\documentclass[11pt,a4paper]{amsart}
\usepackage{etex}

\usepackage{amsmath}
\usepackage{amssymb}
\usepackage{amsthm}
\usepackage{latexsym}
\usepackage{pstricks}
\usepackage{amsbsy}
\usepackage[all]{xypic}
\usepackage{amscd}
\usepackage{mathrsfs}
\usepackage{tipa}
\usepackage{txfonts}
\usepackage{amsfonts}
\usepackage{graphicx}
\usepackage{listings}
\usepackage{url}
\usepackage{bookmark}
\usepackage{tikz}

\newtheorem{thm}{Theorem}[section]
\newtheorem{cor}[thm]{Corollary}
\newtheorem{lem}[thm]{Lemma}
\newtheorem{exm}{Example}
\newtheorem{prop}[thm]{Proposition}
\newtheorem{defn}[thm]{Definition}
\newtheorem{rem}[thm]{Remark}

\newtheorem*{prop*}{Proposition}
\newtheorem*{thm*}{Theorem}

\newenvironment{thmbis}[1]
  {\renewcommand{\thethm}{\ref{#1}$'$}%
   \addtocounter{thm}{-1}%
   \begin{thm}}
  {\end{thm}}

\def\s{\stackrel}

\def\C{\mathcal{C}}
\def\<{\langle}
\def\>{\rangle}

\def\D{\mathcal{D}}

\def\U{\mathscr{U}}

\def\W{\mathscr{W}}
\def\X{\mathscr{X}}
\def\Y{\mathscr{Y}}
\def\Z{\mathcal {Z}}
\def\M{\mathcal {M}}

\def\UU{\mathfrak{U}}
\def\DD{\mathfrak{D}}

\renewcommand{\k}{\mathbf{k}}

\providecommand{\Hom}{\mathop{\rm Hom}\nolimits}%
\providecommand{\Ext}{\mathop{\rm Ext}\nolimits}%
\def\ind{\mbox{ind}}
\def \text{\mbox}

\begin{document}

\title{Mutation of torsion pairs in triangulated categories and its geometric realization}

\author{Yu Zhou}
\address{Department of Mathematical Sciences, Tsinghua University, 100084 Beijing, China}
\curraddr{Department of Mathematical Sciences, Norwegian University of Science and Technology, 7491 Trondheim, Norway}
\email{yu.zhou@math.ntnu.no}

\author{Bin Zhu}
\address{Department of Mathematical Sciences, Tsinghua University, 100084, Beijing, China}
\email{bzhu@math.tsinghua.edu.cn}
\thanks{The authors are supported by the NSF of China (Grants 10771112) and in part by Doctoral Program Foundation of Institute of Higher Education (2009).}

\begin{abstract}
We introduce and study mutation of torsion pairs, as a generalization of mutation of cluster tilting objects, rigid objects and maximal rigid objects. It is proved that any mutation of a torsion pair is again a torsion pair. A geometric realization of mutation of torsion pairs in the cluster category of type $A_n$ or $A_{\infty}$ is given via rotation of Ptolemy diagrams.
\end{abstract}

\maketitle

\section*{Introduction}

Cluster algebras were introduced by Fomin and Zelevinsky \cite{FZ} in order to give an algebraic and combinatorial framework for the positivity and canonical basis of quantum groups. The notion of mutation of clusters is a central idea in the construction of cluster algebras. A categorification of mutation of clusters was introduced by Buan, Marsh, Reineke, Reiten and Todorov \cite{BMRRT}, using mutation of cluster tilting objects in cluster categories. See the nice surveys \cite{Ke2,Ke3,R}. After their work, as generalizations, on the one hand, mutation of cluster tilting objects in arbitrary Hom-finite Krull-Schmidt $\k$-linear triangulated categories was studied in \cite{BIRS,IY,P}. On the other hand, mutation of maximal rigid objects and mutation of rigid objects were studied in \cite{BMV,ZZ} and in \cite{MP}, respectively.

A torsion pair in an abelian category was introduced by Dickson \cite{D}. It is fundamental in algebra, geometry and topology (see the introduction of \cite{BR}) and plays an important role in representation theory of algebras, in particular in tilting theory (see \cite{ASS}). A triangulated version of torsion pairs was studied by Iyama and Yoshino \cite{IY} and a parallel notion, called cotorsion pairs, was studied by Nakaoka \cite{Na} later. Cluster tilting objects, maximal rigid objects and rigid objects are left parts of certain (co)torsion pairs in triangulated categories. 
In the paper, we prove that any mutation of a torsion pair is again a torsion pair. As special cases, we recover mutation of $t$-structures, rigid objects, maximal rigid objects and cluster tilting objects.

Recently, Holm, J{\o}rgensen and M.~Rubey \cite{HJR} gave a bijection from Ptolemy diagrams of a $(n+3)$-gon to torsion pairs in the cluster category of type $A_n$. Later Ng \cite{Ng} gave an analogue bijection from Ptolemy diagrams of an infinity-gon to torsion pairs in the cluster category of type $A_\infty$. In the paper, we define rotation of Ptolemy diagrams and prove that it becomes mutation of torsion pairs under the above bijections.

Since this work, many developments on this topic have been done, e.g. a geometric model of mutation of torsion pairs in cluster categories from unpunctured marked surfaces was given in \cite{ZZZ}, a classification of torsion pairs in 2-Calabi-Yau triangulated categories with cluster tilting objects is given in \cite{ZZ2} and mutation of torsion pairs in cluster category of type $D_n$ was studied in  \cite{G}.

This paper is organized as follows. In section~1, some basic definitions and results concerning torsion pairs are recalled. We introduce the core of a torsion pair and use it to characterize some special cases. In Section~2, we recall from \cite{IY} the construction of subfactor triangulated categories and study its compatibility with (co)torsion pairs. In Section~4, we prove that any mutation of torsion pairs is again a torsion pair and show some properties of mutation of torsion pairs. In the last section, we define rotation of Ptolemy diagrams and prove that it gives a geometric realization of mutation of torsion pairs.

\section{Torsion pairs in triangulated categories}

Let $\k$ be a field and let $\C$ be a triangulated category with the shift functor $[1]$. We denote by $\Hom_{\C}(X,Y)$ (or simply $\Hom(X,Y)$ when no confusion can arise) the set of morphisms from $X$ to $Y$. Through the paper, we always assume that $\C$ satisfies the following conditions:
\begin{enumerate}
  \item[1.] $\C$ is $\k$-linear, i.e. $\Hom(X,Y)$ is a vector space over $\k$ for any objects $X,Y\in\C$, and the compositions of morphisms are $\k$-linear;
  \item[2.] $\C$ is Hom-finite, i.e, $\dim_\k\Hom(X,Y)<\infty$ for any objects $X,Y\in\C$;
  \item[3.] $\C$ is Krull-Schmidt, i.e. every object is isomorphic to a finite direct sum of objects whose endomorphism rings are local.
\end{enumerate}
When we say $\X$ is a subcategory of $\C$, denoted by $\X\subset\C$, we always mean that $\X$ is a full subcategory closed under isomorphisms, direct summands and finite direct sums. For a subcategory $\X$ of $\C$, a \emph{right $\X$-approximation} is a morphism $f: X\rightarrow Y$ in $\C$ with $X\in\X$ such that any morphism from $X'\in \X$ to $Y$ factors through $f$. Such $f$ is called a \emph{right minimal $\X$-approximation} if in addition any endomorphism $g$ of $X$ satisfying $gf=f$ is an isomorphism. We call $\X$ \emph{contravariantly finite} in $\C$, if any object $C\in\C$ admits a right $\X$-approximation. \emph{Left (minimal) $\X$-approximations} and the \emph{covariantly finiteness} of $\X$ can be defined dually. A subcategory $\X$ is called \emph{functorially finite} in $\C$ if $\X$ is both covariantly finite and contravariantly finite in $\C$.

We introduce some notation which will be used later. Let $\X,\Y$ be two subcategories of $\C$. We put
\[\X^{\bot}:=\{M\in\C\mid \Hom(\X,M)=0\}\]
and
\[{}^{\bot}\X:=\{N\in\C\mid \Hom(N,\X)=0\}.\]
We denote by $\X*\Y$ the collection of objects in $\C$ consisting of all such $Z\in\C$ with triangles $X\rightarrow Z\rightarrow Y\rightarrow X[1]$ where $X\in \X$ and $Y\in \Y.$ A subcategory $\X$ of $\C$ is called closed under extensions if $\X\ast\X=\X$.

\begin{defn}
Let $\X$, $\Y$ be subcategories of $\C$.
\begin{itemize}
  \item[(1)]\cite[Definition~2.2]{IY} The pair $(\X,\Y)$ is called a torsion pair if
  \[\Hom(\X,\Y)=0\ \text{and}\ \C=\X*\Y.\]
  \item[(2)]\cite[Definition~2.1]{Na} The pair $(\X,\Y)$ is called a cotorsion pair if
  \[\Hom(\X,\Y[1])=0\ \text{and}\ \C=\X*\Y[1].\]
\end{itemize}
\end{defn}

It is obvious that a pair $(\X,\Y)$ of subcategories is a cotorsion pair if and only if $(\X, \Y [1])$ is a torsion pair. So there is a bijection between torsion pairs and cotorsion pairs, and a subcategory $\X$ is the left part of a torsion pair if and only if it is the left part of a cotorsion pair. Hence, to study torsion pairs is equivalent to studying cotorsion pairs.

\begin{prop}[\cite{IY}, Proposition~2.3]\label{prop:iy}
If $(\X,\Y)$ is a cotorsion pair, then $\X={}^{\bot}(\Y[1])$ and $\Y=(\X[-1])^\bot$. For a subcategory $\X$ of $\C$, the pair $(\X, (\X[-1])^\bot)$ is a torsion pair if and only if $\X$ is contravariantly finite in $\C$ and is closed under extensions.
\end{prop}


\begin{defn}
Let $\X$ be a subcategory of $\C$.
\begin{itemize}
  \item[(1)] We call $\X$ rigid if $\Hom(\X,\X[1])=0$.
  \item[(2)] We call $\X$ maximal rigid provided that $\X$ is rigid and any rigid subcategory of $\C$ containing $\X$ is $\X$.
  \item[(3)] We call $\X$ cluster tilting if $\X$ is contravariantly finite, and for any $X\in\C$, we have that $X\in \X$ if and only if $\Hom(\X,X[1])=0$.
\end{itemize}
\end{defn}

Functorially finite rigid subcategories, functorially finite maximal rigid subcategories and cluster tilting subcategories are received much attention recently \cite{BMRRT,BIRS,BMV,BR,IY,KR,KZ,MP,ZZ}.

\begin{rem}
The original definition of cluster tilting subcategories in \cite{KR} requires more conditions than here. It turns out by \cite[Lemma~3.2]{KZ} that the present one is equivalent to the original one. In particular, cluster tilting subcategories are always functorially finite and rigid.
\end{rem}

We introduce some related special classes of torsion pairs.

\begin{defn}\label{def:cases}
Let $(\X,\Y)$ be a (co)torsion pair in $\C$.
\begin{itemize}
  \item [(1)] We call $(\X,\Y)$ a (co)torsion t-structure if $\X[1]\subset\X$.
  \item [(2)] We call $(\X,\Y)$ rigid (or maximal rigid, or cluster tilting, respectively) if $\X$ is rigid (or maximal rigid, or cluster tilting, respectively).
\end{itemize}
\end{defn}

\begin{rem}\label{rmk:cases}
By definition, we have the following facts.
\begin{itemize}
  \item [(1)] A pair $(\X,\Y)$ of subcategories of $\C$ is a torsion (or cotorsion, respectively) t-structure if and only if $(\X,\Y[1])$ (or $(\X,\Y[2])$, respectively) is a t-structure in the sense of \cite{BBD}. So there is a bijection between (co)torsion t-structures and t-structures, and to study t-structures is equivalent to studying (co)torsion t-structures.
  \item [(2)] Any functorially finite rigid (or maximal rigid, or cluster tilting, respectively) subcategories are contravariantly finite and closed under extensions. Then by Proposition~\ref{prop:iy} they are the left part of rigid (or maximal rigid, or cluster tilting, respectively) (co)torsion pairs.
\end{itemize}
\end{rem}

Let $\X$ be a subcategory of $\C$. An object $X\in\X$ is called $\Ext$-injective (resp. $\Ext$-projective) in $\X$ provided that $\Hom(\X, X[1])=0$ (resp. $\Hom(X, \X[1])=0$). The subcategory of $\X$ consisting of all the $\Ext$-injective (resp. $\Ext$-projective) objects in $\X$ is denoted by $I(\X)$ (resp. $P(\X)$). We introduce the following notion.

\begin{defn}
Let $(\X,\Y)$ be a (co)torsion pair. We call $I(\X)$ the \emph{core} of $(\X,\Y)$.
\end{defn}

The core of a (co)torsion pair has nice properties.

\begin{prop}\label{prop:core}
Let $(\X,\Y)$ be a cotorsion pair of $\C$. Then $I(\X)$ is rigid and $I(\X)=\X\cap \Y=P(\Y)$. 
Moreover, for any $X\in\X$ and $Y\in\Y$, there are triangles
\[X\xrightarrow{g} E\to X'\to X[1] \quad \text{and} \quad Y'\to F\xrightarrow{h} Y\to Y'[1]\]
with $g$ a left $I(\X)$-approximation, $h$ a right $I(\X)$-approximation, $X'\in\X$ and $Y'\in\Y$.
\end{prop}

\begin{proof}
Since $I(\X)$ is Ext-injective in $\X$, We have $\Hom(\X, I(\X)[1])=0$. In particular, $\Hom(I(\X),I(\X)=0$, so $I(\X)$ is rigid. Let $X\in \X$. Then $\Hom(\X, X[1])=0$ if and only if $X\in \Y$ by Proposition~\ref{prop:iy}. Then $X\in I(\X)$ if and only if $X\in \X\cap \Y$. Hence $I(\X)=\X\cap \Y$. Similarly, we have $P(\Y)=\X\cap\Y$.

Since $(\X,\Y)$ is a cotorsion pair, we have that $(\X[-1],\Y[-1])$ is also a cotorsion pair. Then for any $X\in \X$, there exists a triangle $X_1[-1]\rightarrow X\s{g}{\rightarrow} Y_1\rightarrow X_1$ with $X_1\in \X$ and $Y_1\in \Y$. It follows that $Y_1\in \X$ since $\X$ is closed under extensions. Hence $Y_1\in I(\X)$. So $g$ is a left $I(\X)$-approximation of $X$. Thus we get the first required triangle. Similarly, we have the second triangle.
\end{proof}

The special cases of (co)torsion pairs in Definition~\ref{def:cases} can be characterized by their cores.

\begin{prop}\label{prop:cases}
Let $(\X,\Y)$ be a cotorsion pair in $\C$. Then the following hold.
\begin{enumerate}
\item $(\X,\Y)$ is a cotorsion $t$-structure if and only if $I(\X)=0.$
\item $(\X,\Y)$ is rigid if and only if $I(\X)=\X$.
\item $(\X,\Y)$ is cluster tilting if and only if $I(\X)=\X=\Y$.
\end{enumerate}
\end{prop}

\begin{proof}
\begin{enumerate}
\item By definition, it is sufficient to prove that $\X[1]\subset\X$ if and only if $I(\X)=0$. If $\X[1]\subset\X$, we have that $I(\X)=\X\cap \Y\subseteq\X[-1]\cap \Y$. Since $\Hom(\X[-1],\Y)\cong\Hom(\X,\Y[1])=0$, $\X[-1]\cap\Y$ have to be zero. So $I(\X)=0$. Conversely, suppose that $I(\X)=0$. For any $X\in\X$, by Proposition~\ref{prop:core}, there is a triangle $X\to E\to X'\to X[1]$ with $E\in I(\X)$ and $X'\in\X$. So $E=0$ and hence $X[1]\cong X'\in\X$. This proves that $\X[1]\subseteq \X$.
\item If $\X$ is rigid, then any object in $\X$ is Ext-injective. So it follows that $I(\X)=\X$. Conversely, if $I(\X)=\X$, then $\X$ is rigid by Proposition~\ref{prop:core}.
\item By definition, $\X$ is cluster tilting if and only if $\X=\{X\in\C\mid \Hom(\X,X[1])=0\}$. On the other hand, we have that $\Y=\{X\in\C\mid \Hom(\X,X[1])=0\}$ by Proposition~\ref{prop:iy}. Hence $\X$ is cluster tilting if and only if $\X=\Y$. In this case, $I(\X)=\X=\Y$.
\end{enumerate}
\end{proof}

\section{Compatibility with subfactor triangulated categories}

From now on, we assume that $\C$ has a Serre functor $\mathbb{S}$. Recall from \cite{BK} that a functor $\mathbb{S}:\C\rightarrow \C$ is called a Serre functor if there is a bifunctorial isomorphism
\[\Hom(X,Y)\cong D \Hom(Y,\mathbb{S} X)\]
for any objects $X,Y\in \C$, where $D=\Hom_\k(-,\k)$. 

In this section, we fix a functorially finite subcategory $\D$ of $\C$ which satisfying the condition $\mathbb{S}\D[-2]=\D$. Then it is easy to check that ${}^\bot\D[1]=\D[-1]^\bot$, which is denoted by $\Z$. The quotient category $\U:=\Z/\D$, called a subfactor triangulated category, is defined by the following data: $\U$ has the same objects as $\Z$, and for each $X,Y\in\Z$, the morphism space $\Hom_\U(X,Y)$ is defined as
\[\Hom_\U(X,Y):=\Hom_\Z(X,Y)/\D(X,Y),\]
where $\D(X,Y)$ is the subspace of $\Hom_\Z(X,Y)$ consisting of morphisms factoring through objects in $\D$. It is proved in \cite{IY} that $\U$ carries a natural triangulated structure inherited from the triangulated structure of $\C$ as follows:
\begin{itemize}
  \item for any object $X$, choose a left $\D$-approximation $f: X\rightarrow D$ and extend it to a triangle \[X\s{f}{\longrightarrow}D\longrightarrow Z\longrightarrow X[1],\]
      then the shift of $X$ in $\U$ is defined to be $Z$, denoted by $X\<1\>$;
  \item for any triangle $X\s{a}{\rightarrow} Y\s{b}{\rightarrow}
Z\s{c}{\rightarrow} X[1]$ in $\C$ with $X,Y,Z\in\Z$, we have a commutative diagram:
\[\xymatrix{
X\ar[r]^a\ar@{=}[d]&Y\ar[r]^b\ar[d]&Z\ar[r]^c\ar[d]^d&X[1]\ar@{=}[d]\\
X\ar[r]^f&D\ar[r]^g&X\langle1\rangle\ar[r]^h&X[1]
}\]
where the existence of $d$ follows from $\Hom_\C(Z,D[1])=0$. Considering the complex
\[X\xrightarrow{\overline{a}}Y\xrightarrow{\overline{b}}Z\xrightarrow{\overline{d}}X\langle1\rangle\]
where $\overline{a}, \overline{b}, \overline{d}$ are the images, under the quotient functor $\Z\rightarrow \U$, of the maps $a, b, d$ respectively, the triangles in $\U$ are defined to be the complexes (up to isomorphism) obtained by this way.
\end{itemize}

The following lemma will be useful later.



\begin{lem}\label{lem:eq}
Let $(\X,\Y)$ be a cotorsion pair. Then $\D\subset I(\X)$ if and only if $\D\subset\X\subset\Z$ and $\D\subset\Y\subset\Z$.
\end{lem}

\begin{proof}
The 'if' part is trivial by $I(\X)=\X\cap\Y$. To prove the 'only if' part, let $X$ be an arbitrary object in $\X$. Since $\D\subset I(\X)\subset \Y$, we have $\Hom(X,\D[1])=0$. It follows that $X\in\Z$. Hence $\X\subset\Z$. Similarly, we have $\Y\in\Z$. The proof is complete.
\end{proof}

Thus, one can consider the relationship between cotorsion pairs in $\C$ whose cores contain $\D$ and cotorsion pairs in the subfactor triangulated category $\U$. In the following, $\overline{\W}$ denotes the subcategory of $\U$ consisting of objects $W\in \W$, for a subcategory $\W$ satisfying $\D\subset\W\subset\Z$. It is clear that any subcategory of $\U$ has this form.

\def\bX{\overline{\X}}
\def\bY{\overline{\Y}}

\begin{thm}\label{thm:ud}
Let $(\X,\Y)$ be a cotorsion pair in $\C$ with $\D\subset I(\X)$. Then $(\overline{\X},\overline{\Y})$ is a cotorsion pair in $\U$ with $I(\overline{\X})=\overline{I(\X)}$. Moreover, the map $(\X,\Y)\mapsto(\overline{\X},\overline{\Y})$ is a bijection from the set of cotorsion pairs in $\C$ whose cores contain $\D$ and the set of cotorsion pairs in $\U$.
\end{thm}

\begin{proof}
Since $(\X,\Y)$ is a cotorsion pair with $\D\subset I(\X)$, we have $\D\subset\X\subset\Z$ and $\D\subset\Y\subset\Z$ by Lemma~\ref{lem:eq}. So there are subcategories $\overline{\X}$ and $\overline{\Y}$ of $\U$. By \cite[Lemma~4.8]{IY}, we have an isomorphism $\Hom_\C(\X,\Y[1])\cong\Hom_\U(\bX,\bY\<1\>)$. It follows that $\Hom_\U(\bX,\bY\<1\>)=0$. Since $\X\ast\Y[1]=\C$, for any $Z\in\Z$, there is a triangle
$Y\rightarrow X\rightarrow Z\rightarrow Y[1]$ in $\C$
with $X\in\X$ and $Y\in\Y$. Since all of $Y$, $X$, $Z$ are in $\Z$, by the construction of triangles in $\U$, there is a triangle
$Y\rightarrow X\rightarrow Z\rightarrow Y\langle1\rangle$ in $\U$. Therefore, we have $\U=\overline{\X}\ast\overline{\Y}\langle1\rangle$ and then $(\bX,\bY)$ is a cotorsion pair in $\U$. In this case, we have $I(\overline{\X})=\overline{\X}\cap\overline{\Y}=\overline{\X\cap\Y}=\overline{I(\X)}$.

Let $(\X_i,\Y_i)$ be cotorsion pairs in $\C$ with $I(\X_i)\subset\D$, for $i=1,2$. If $\bX_1=\bX_2$, then by definition, $\X_1=\X_2$. Hence the map $(\X,\Y)\mapsto(\overline{\X},\overline{\Y})$ is injective. Now consider subcategories $\X$ and $\Y$ of $\C$ satisfying $\D\subset\X\subset\Z$ and $\D\subset\Y\subset\Z$ such that $(\bX,\bY)$ is a cotorsion pair in $\U$. By definition, we have $\Hom_\U(\bX,\bY)=0$ and $\U=\bX\ast\bY\<1\>$. Using \cite[Lemma~4.8]{IY} again, we have $\Hom_\C(\X,\Y[1])=0$. So by \cite[Proposition~2.1]{IY}, $\X\ast\Y[1]$ is closed under direct summands. On the other hand, for any
$Z\in\Z$, there is a triangle $Y\rightarrow X\rightarrow Z\rightarrow Y\langle1\rangle$ in $\U$
with $X\in\bX$ and $Y\in\bY$. Then by the construction of triangles in $\U$, there is a triangle $Y\rightarrow X\rightarrow Z'\rightarrow Y[1]$
in $\C$ such that $Z\cong Z'$ in $\U$. Then $Z\cong Z'$ in $\C$ up to direct summands in $\D$. Since $\X\ast\Y[1]$ is closed under direct summands, we have $Z\in\X\ast\Y[1]$ and hence $\Z\subset\X\ast\Y[1]$. Since $\D$ is functorially finite in $\C$ and is closed under extensions, $(\Z,\D)$ is a cotorsion pair by Proposition~\ref{prop:iy}. Hence $\C=\Z\ast\D[1]\subset\X\ast\Y[1]\ast\D[1]=\X\ast\Y[1]$. Therefore, $(\X,\Y)$ is a cotorsion pair in $\C$. By Lemma~\ref{lem:eq}, $\D\subset I(\X)$. It follows that the map $(\X,\Y)\mapsto(\overline{\X},\overline{\Y})$ is surjective and then it is a bijection.
\end{proof}

Restricting to the special cases of cotorsion pairs in Definition~\ref{def:cases}, we have the following bijections where the last one is known by \cite[Theorem~5.1]{IY}.

\begin{cor}\label{cor:cases}
The bijection in Theorem~\ref{thm:ud} induces the following bijections:
\begin{enumerate}
\item a bijection between cotorsion pairs in $\C$ whose cores are $\D$ and cotorsion $t$-structure in $\U$,
\item a bijection between rigid cotorsion pairs in $\C$ whose cores contains $\D$ and rigid cotorsion pairs in $\U$,
\item a bijection between maximal rigid cotorsion pairs in $\C$ whose cores contains $\D$ and maximal rigid cotorsion pairs in $\U$, and
\item a bijection between cluster tilting cotorsion pairs in $\C$ whose cores contains $\D$ and cluster tilting cotorsion pairs in $\U$.
\end{enumerate}
\end{cor}

\begin{proof}
It follows directly from Proposition~\ref{prop:cases} and Theorem~\ref{thm:ud} that we have the bijections in (1), (2), and (4). Now we prove the bijection in (3). By Theorem~\ref{thm:ud}, we only need to prove that $\X$ is maximal rigid if and only if $\bX$ is maximal rigid. Assume that $\X$ is maximal rigid. For any rigid subcategory $\X'$ of $\Z$ with $\bX\subset\overline{\X'}$, we have that $\X\subset\X'$. It follows from the bijection in (2) that $\X'$ is rigid in $\C$. So $\X=\X'$ and then $\bX=\overline{\X'}$ which implies that $\bX$ is maximal rigid. Conversely, assume that $\bX$ is maximal rigid. For any rigid subcategory $\X'$ of $\C$ with $\X\subset\X'$, we have $\D\subset \X'$. Since $\X'$ is rigid, we have $\Hom_\C(\X',\X'[1])=0$. In particular, we have $\Hom_\C(\X',\D[1])=0$. So $\X'\in{}^\bot\D[1]=\Z$. Then $\bX\subset\overline{\X'}$ with $\overline{\X'}$ rigid in $\U$. So $\bX=\overline{\X'}$, and then $\X=\X'$ which implies that $\X$ is maximal rigid.

\end{proof}

\section{Mutation of torsion pairs}

In this section, we study mutation of (co)torsion pairs. First we recall from \cite{IY} the notion of mutation of subcategories. Let $\C$ be a Krulll-Schmidt Hom-finite $\k$-linear triangulated category with a Serre functor $\mathbb{S}$.

\begin{defn}[\cite{IY}, Definition~2.5]
Fix a functorially finite rigid subcategory $\D$ of $\C$. For any object $M\in\C$, the forward and backward $\D$-mutations $\mu^{-1}_\D(M)$ and $\mu_\D(M)$ are defined respectively by the following triangles
\[M\xrightarrow{f}D'\rightarrow \mu^{-1}_\D(M)\rightarrow M[1]\]
and
\[M[-1]\rightarrow \mu_\D(M)\rightarrow D''\xrightarrow{g} M\]
with $f$ a left minimal $\D$-approximation and with $g$ a right minimal $\D$-approximation. Furthermore, for any subcategory $\M$ of $\C$ with $\D\subset\M$, the forward and backward $\D$-mutations $\mu^{-1}_\D(\M)$ and $\mu_\D(\M)$ are defined respectively by
\[\mu^{-1}(\M;\D):=\{\mu^{-1}_{\D}(M)\mid M\in \M\}\cup\D\]
and
\[\mu(\M;\D):=\{\mu_{\D}(M)\mid M\in \M\}\cup\D.\]
\end{defn}

It is clear that when $\D=0$, we have $\mu^{-1}(\M;0)=\M[1]$ and $\mu(\M;0)=\M[-1]$.

\begin{prop}[\cite{IY}, Proposition~2.7]\label{rmk:mu}
Let $\D$ be a functorially finite rigid subcategory of $\C$ with $\mathbb{S}\D[-2]=\D$ and let $\Z={}^\bot\D[1]$. Then the maps $\mu^{-1}(-;\D)$ and $\mu(-;\D)$ are mutually inverse on the set of subcategories $\M$ of $\C$ satisfying $\D\subset\M\subset\Z$.
\end{prop}





The follow theorem is a generalization of well-known results on mutation of cluster tilting objects \cite{BMRRT,IY}, rigid objects \cite{MP} and maximal rigid objects \cite{BIRS,BMV,ZZ} to the setting of cotorsion pairs.

\begin{thm}\label{thm:mutation}
Let $(\X,\Y)$ be a cotorsion pair in $\C$ and let $\D\subset I(\X)$ be a functorially finite subcategory of $\C$ satisfying $\mathbb{S}\D=\D[2]$. Then the pairs $(\mu^{-1}(\X;\D),\mu^{-1}(\Y;\D))$ and $(\mu(\X;\D),\mu(\Y;\D))$ are cotorsion pairs in $\C$ and we have $I(\mu^{-1}(\X;\D))=\mu^{-1}(I(\X);\D)$ and $I(\mu(\X;\D))=\mu(I(\X);\D)$.
\end{thm}

\begin{proof}
By Proposition~\ref{prop:core}, we have that $I(\X)$ is rigid, then so is $\D$. 
It follows that both $\mu^{-1}(\X;\D)$ and $\mu^{-1}(\Y;\D)$ are defined. 
We denote $\mu^{-1}(\X;\D)$, $\mu^{-1}(\Y;\D)$ and $\mu^{-1}(I(\X);\D)$ 
by $\X'$, $\Y'$ and $I'$ respectively. 
As in the previous section, 
we denote by $\U$ the subfactor triangulated category $\Z/\D$, 
where $\Z=^\bot\D[1]$. Then by Theorem~\ref{thm:ud}, 
we have that $(\overline{\X},\overline{\Y})$ is a cotorsion pair in $\U$ 
and $I(\bX)=\overline{I(\X)}$. Note that its $0$-mutation 
$(\overline{\X}\langle1\rangle,\overline{\Y}\langle1\rangle)$ in $\U$ is also a cotorsion pair. 
By \cite[Proposition 4.4]{IY}, we have $\overline{\X}\langle1\rangle=\overline{\X'}$ and $\overline{\Y}\langle1\rangle=\overline{\Y'}$. Then by Theorem~\ref{thm:ud} again, we have that $(\X',\Y')$ is a cotorsion pair in $\C$ and $I(\overline{\X'})=\overline{I(\X')}$. So
$\overline{I'}=\overline{I(\X)}\langle1\rangle=I(\overline{\X})\langle1\rangle=I(\overline{\X}\langle1\rangle)=I(\overline{\X'})=\overline{I(\X')}.$
Therefore, $I(\mu^{-1}(\X;\D))=\mu^{-1}(I(\X);\D)$. The assertion for $(\mu(\X;\D),\mu(\Y;\D))$ can be proved similarly.
\end{proof}

Recall that $(\X,\Y)$ is a torsion pair if and only if $(\X,\Y[-1])$ is a cotorsion pair, and note that $\mu^{-1}(\Y;\D[1])[-1]=\mu^{-1}(\Y[-1];\D)$ and $\mu(\Y;\D[1])[-1]=\mu^{-1}(\Y[-1];\D)$. So we have an analogue theorem for torsion pairs.

\begin{thmbis}{thm:mutation}\label{thm:mutation2}
Let $(\X,\Y)$ be a torsion pair in $\C$ and 
let $\D\subset I(\X)$ be a functorially finite subcategory of $\C$ satisfying $\mathbb{S}\D=\D[2]$. 
Then the pairs $(\mu^{-1}(\X;\D),\mu^{-1}(\Y;\D[1]))$ and $(\mu(\X;\D),\mu(\Y;\D[1]))$ 
are torsion pairs in $\C$ and we have $I(\mu^{-1}(\X;\D))=\mu^{-1}(I(\X);\D)$ 
and $I(\mu(\X;\D))=\mu(I(\X);\D)$.
\end{thmbis}

Now we introduce the notion of mutation of (co)torsion pairs.

\begin{defn}\label{defn:mut}
Let $(\X,\Y)$ be a (co)torsion pair in $\C$ and let $\D\subset I(\X)$ be a functorially finite subcategory of $\C$ with $\mathbb{S}\D[-2]=\D$.
\begin{itemize}
  \item[(1)] When $(\X,\Y)$ is a torsion pair, the torsion pair $(\mu ^{-1}(\X;\D),\mu^{-1}(\Y;\D[1]))$ is called the forward $\D$-mutation of $(\X,\Y)$ and the torsion pair $(\mu(\X;\D),\mu(\Y;\D[1]))$ is called the backward $\D$-mutation of $(\X,\Y)$.
  \item[(2)] When $(\X,\Y)$ is a cotorsion pair, the cotorsion pair $(\mu ^{-1}(\X;\D),\mu^{-1}(\Y;\D))$ is called the forward $\D$-mutation of $(\X,\Y)$ and 
  the cotorsion pair $(\mu(\X;\D),\mu(\Y;\D))$ is called the backward $\D$-mutation of $(\X,\Y)$
\end{itemize}
\end{defn}

It follows from Proposition~\ref{rmk:mu} that if $(\X,\Y), (\X',\Y')$ are two (co)torsion pairs in $\C$, then $(\X',\Y')$ is the forward $\D$-mutation of $(\X,\Y)$ if and only if $(\X,\Y)$ is the backward $\D$-mutation of $(\X',\Y')$. So we can define the following equivalent relation.

\def\mut{\rm Mut}

\begin{defn}\label{defn:equiv}
Two (co)torsion pairs in $\C$ are called mutation equivalent if we can get one from the other one by a sequence of forward or backward mutations. We denote by $\mut((\X,\Y))$ the equivalence class containing a (co)torsion pair $(\X,\Y)$ under this mutation equivalent relation.
\end{defn}

Now we apply mutation of torsion pairs to the special cases in Definition~\ref{def:cases}. Mutations of such special cases appeared recently in the study of cluster tilting theory and cluster algebras, see the surveys \cite{Ke2,R} and the references there.

\begin{cor}\label{cor:cases2}
Let $(\X,\Y)$ be a (co)torsion pair. If $(\X,\Y)$ is a (co)torsion t-structure, (or a rigid (co)torsion pair, or a maximal rigid (co)torsion pair, or a cluster tilting (co)torsion pair, respectively), then so is any (co)torsion pair in $\mut((\X,\Y))$.
\end{cor}

\begin{proof}
We only consider the case that $(\X,\Y)$ is a cotorsion pair. Let $(\X',\Y')$ and $(\X'',\Y'')$ be the forward $\D$-mutation and the backward $\D$-mutation of $(\X,\Y)$ respectively for some $\D$. It is sufficient to prove the assertion for $(\X',\Y')$ and $(\X'',\Y'')$. By \cite[Proposition~4.4]{IY}, we have $(\overline{\X'},\overline{\Y'})=(\overline{\X}\langle1\rangle,\overline{\Y}\langle1\rangle)$ in the subfactor triangulated category $\U$. Note that the shift of a cotorsion t-structure (or rigid cotorsion pair, or maximal rigid cotorsion pair, or cluster tilting cotorsion pair, respectively) is again a cotorsion t-structure (or rigid cotorsion pair, or maximal rigid cotorsion pair, or cluster tilting cotorsion pair, respectively). Hence the assertion for $(\X',\Y')$ follows from Corollary~\ref{cor:cases}. The proof of the assertion for $(\X'',\Y'')$ is similar.
\end{proof}

In the following, we study properties of mutation of (co)torsion pairs. Let $\ind\X$ denote the set of isoclasses of indecomposable objects in a subcategory $\X$ of $\C$.

\begin{prop}\label{prop:prop}
Let $(\X',\Y')$ be a forward $\D$-mutation of a (co)torsion pair $(\X,\Y)$. Then the following hold.
\begin{enumerate}
  \item There is a bijection between $\ind\X$ and $\ind\X'$, which restricting to $I(\X)$ induces a bijection between $\ind I(\X)$ and $\ind I(\X')$.
  \item If $\D$ is a proper subcategory of $I(\X)$ (i.e. $\D\neq I(\X)$), then $I(\X')\neq I(\X)$ and $(\X',\Y')\neq(\X,\Y)$.
  \item If $\D=I(\X)$, then we have $I(\X')=I(\X)$, $\X'\subseteq \X$ and $\Y\subseteq \Y'$.
\end{enumerate}
\end{prop}

\begin{proof}
We still only consider the case that $(\X,\Y)$ is a cotorsion pair.
\begin{enumerate}
  \item Since $(\X',\Y')$ is a forward $\D$-mutation of $(\X,\Y)$, 
  we have $\overline{\X'}=\overline{\X}\langle1\rangle$ and 
  $\overline{I(\X')}=\overline{I(\X)}\langle1\rangle$ in $\U$. 
  Then the shift functor $\<1\>$ induces a bijection between 
  $\ind\overline{\X}$ and $\ind\overline{\X'}$, 
  which restricting to $\overline{I(\X)}$ induces a bijection between $\ind\overline{I(\X)}$ and $\ind\overline{I(\X')}$.  Note that $\ind\W=\ind\overline{\W}\cup\ind\D$ for any subcategory $\W$ of $\C$ with $\D\subset\W\subset\Z$. Thus these bijections can be lifted to the required bijections.
  \item We prove that $I(\X')\neq I(\X)$, then $(\X',\Y')\neq(\X,\Y)$ follows. Since $\D$ is a proper subcategory of $I(\X)$, there exists an object $X\in I(\X)$ which does not belong to $\D$. Let $f: X\rightarrow D$ be a minimal left $\D$-approximation of $X$. Then we have a triangle
      \[X\xrightarrow{f} D\to \mu^{-1}_\D(X)\to X[1].\]
      By definition, $\mu^{-1}_\D(X)$ is an object in $\mu^{-1}(I(\X);\D)$. So by Theorem~\ref{thm:mutation}, we have $\mu^{-1}_\D(X)\in I(\X')$. We claim that $\mu^{-1}_\D(X)\notin I(\X)$, and then $I(\X)\neq I(\X')$. Indeed, if $\mu^{-1}_\D(X)\in I(\X)$, since $I(\X)$ is rigid by Proposition~\ref{prop:core}, the above triangle is split. So the map $f$ is a section, which is a contradiction to that $f$ is left minimal.
  \item If $\D=I(\X)$, then we have $I(\X')=I(\X)=\D$. So $(\overline{\X},\overline{\Y})$ is a cotorsion $t$-structure in $\U$ by Corollary~\ref{cor:cases}. It follows that $\overline{\X'}=\overline{\X}\langle1\rangle\subseteq\overline{\X}$ which implies that $\X'\subset\X$. Then $\Y=\X[-1]^\bot\subset\X'[-1]^\bot=\Y'$.
\end{enumerate}

\end{proof}

A triangulated category is called of finite type if it has finitely many indecomposable objects up to isomorphism. Triangulated categories of finite type were studied recently in \cite{XZ1,XZ2,Am,Koe,Kr}. Cluster categories of Dynkin type \cite{BMRRT} and stable Cohen-Macaulay categories of finite type \cite{BIKR} provide examples of triangulated categories of finite type.

\medskip

\begin{cor}\label{cor:finite}
Let $(\X',\Y')$ be a forward $\D$-mutation of a (co)torsion pair $(\X,\Y)$. If $\C$ is of finite type, then $(\X',\Y')=(\X,\Y)$ if and only if $I(\X')=I(\X)$.
\end{cor}

\begin{proof} The "only if" part is obvious. We prove the "if" part. Suppose that $I(\X')=I(\X)$, then by Proposition~\ref{prop:prop}(2) we have $\D=I(\X)$. So by Proposition~\ref{prop:prop}(3), $\X'\subseteq\X$. On the other hand, by Proposition~\ref{prop:prop}(1), there is a bijection between $\ind\X$ and $\ind\X'$. So we have $\ind\X=\ind\X'$ since $\C$ is of finite type. Therefore $\X'=\X$, and then $(\X',\Y')=(\X,\Y)$.
\end{proof}

We note that there are different torsion pairs having the same core in general. For example, the trivial torsion pairs $(\C,0)$ and $(0,\C)$ have the same core 0. We also note that Corollary~\ref{cor:finite} does not hold without the condition that $\C$ is of finite type. For example we take $\C$ to be the bounded derived category of category of finitely generated right modules of a finite dimensional $\k$-algebra. Let $(\X,\Y)$ be the canonical $t$-structure in $\C$ and let $\D=0$. Then $(\X,\Y[-1])$ is a torsion $t$-structure and its forward $\D$-mutation $(\X',\Y')$ is $(\X[1], \Y)$. By Proposition~\ref{prop:cases}(1), the cores of these two torsion pairs are 0. But we have $\X'\not=\X$.

\section{A geometric realization of mutation of torsion pairs}\label{sec:geo}

In \cite{HJR} and \cite{Ng}, the torsion pairs in cluster categories of types $A_n$ and $A_\infty$ are classified in terms of Ptolemy diagrams of polygons, respectively. In this section, we introduce the notion of rotation of Ptolemy diagrams, which gives a geometric realization of mutation of torsion pairs in such categories. Note that cluster categories are 2-Calabi-Yau, that is, $[2]$ is a Serre functor. Hence any subcategory $\D$ satisfies the condition $\mathbb{S}\D[-2]=\D$.

\subsection{Dynkin type A}

\def\CA{{\C_{A_n}}}

Let $P_{n+3}$ be a convex polygon with $n+3$ vertices labeled by $1, 2, \cdots, n+3\in \mathbb{Z}/(n+3)\mathbb{Z}$ in counterclockwise order. A line segment in $P_{n+3}$ linking two non-adjacent (resp. adjacent) vertices $i,j$ is called a diagonal (resp. an edge), and we denote it by $\{i,j\}$. For any diagonal $\{i,j\}$, its rotation in $P_{n+3}$ is defined to be $\{i-1,j-1\}$.

We refer to \cite{BMRRT} for the definition of the cluster category $\C_{A_n}$ of type $A_n$. By \cite{CCS}, there is a bijection from the set of diagonals of $P_{n+3}$ to the set $\ind\C_{A_n}$ of isoclasses of indecomposable objects in $\C_{A_n}$ and under this bijection, the rotation of diagonals becomes the shift of objects. We denote by $X_{i,j}$ the object corresponding to a diagonal $\{i,j\}$ (so we have $X_{i,j}=X_{j,i}$). Since in the paper any subcategory is assumed to be full and to be closed under isomorphisms, direct summands and finite direct sums, the bijection above induces a bijection from sets of diagonals of $P_{n+3}$ to subcategories of $\C_{A_n}$. For a set $\UU$ of diagonals, we denote by $X(\UU)$ the corresponding subcategory of $\C_{A_n}$.

\begin{defn}[\cite{HJR}]
Let $\UU$ be a set of diagonals of $P_{n+3}$. It is called a Ptolemy diagram if it satisfies that when there are two crossing diagonals $\{i,j\}$ and $\{k,l\}$ in $\UU$, then those of $\{i,k\}$, $\{k,j\}$, $\{j,l\}$, $\{l,i\}$ which are diagonals are in $\UU$. See Figure~\ref{fig:Pr}.
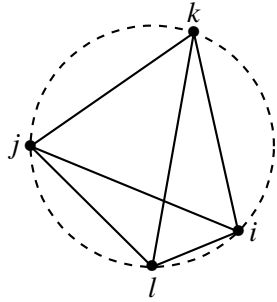
\begin{figure}[ht]\centering
\begin{tikzpicture}[scale=.8]
\draw[thick, dashed](0,0) circle (2);
\draw[] (-2,0)node{$\bullet$} (70:2)node{$\bullet$} (-45:2)node{$\bullet$} (-90:2)node{$\bullet$};
\draw[thick] (-2,0)node[left]{$j$} (70:2)node[above]{$k$} (-45:2)node[right]{$i$} (-90:2)node[below]{$l$};
\draw[thick] (-2,0) to (70:2) to (-45:2) to (-90:2) to (-2,0) to (-45:2);
\draw[thick] (70:2) to (-90:2);
\end{tikzpicture}
\caption{The Ptolemy relation for a polygon  }
\label{fig:Pr}
\end{figure}
\end{defn}

\begin{prop}[\cite{HJR}]\label{prop:hjr}
For a set of diagonals $\UU$, the pair $(X(\UU),X(\UU)^{\bot})$ is a torsion pair in $\C_{A_n}$ if and only if the set $\mathfrak{U}$ is a Ptolemy diagram.
\end{prop}

Now we define frame and mutation of a Ptolemy diagram.

\begin{defn}
Let $\UU$ be a Ptolemy diagram. The set of diagonals in $\UU$ that does not cross any diagonals in $\UU$ is called the frame of $\UU$, denoted by $F_\UU$.
\end{defn}

\begin{defn}
Let $\DD$ be a set of non-crossing diagonals.
\begin{itemize}
  \item[(1)] The diagonals in $\DD$ divide the polygon $P_{n+3}$ into several polygons, called $\DD$-cells. Then any diagonal $\{i,j\}$ which neither is in $\DD$ nor crosses any diagonals in $\DD$ is a diagonal of a $\DD$-cell. We call the rotation of $\{i,j\}$ in this $\DD$-cell the $\DD$-rotation of $\{i,j\}$ and denote it by $\rho_\D(\{i,j\})$.
  \item[(2)] Let $\UU$ be a Ptolemy diagram with $\DD\subset F_\UU$. The $\DD$-rotation of $\UU$ is defined as
  \[\rho(\UU;\DD):=\{\rho_\DD(\{i,j\})\mid \{i,j\}\in\UU\setminus \DD\}\cup \DD.\]
\end{itemize}
\end{defn}

\begin{exm} Let $P$ be a polygon with $8$ vertices, and let \[\mathfrak{U}=\{ \{2,7\},\{2,8\},\{3,7\},\{3,8\},\{4,6\},\{4,7\},\{5,7\}\}\] be a Ptolemy diagram of $P_{n+3}$, see the left picture in Figure~\ref{fig:rP}. Then $F_\UU=\{\{2,8\},\{3,7\},\{4,7\}\}$. Let $\DD=\{\{3,7\}\}\subset F_\UU$. Then $\DD$ divides $P$ into two polygons. The $\D$-mutation of $\UU$ is the right picture in Figure~\ref{fig:rP}, where for example, $\rho_\DD(\{2,7\})=\{1,3\}$. It is clear that $\rho(\UU;\DD)$ is again a Ptolemy diagram.
\begin{figure}[ht]\centering
\begin{tikzpicture}[scale=.8]
\draw[thick] (22.5:2)node{$\bullet$} to (67.5:2)node{$\bullet$} to (112.5:2)node{$\bullet$}  to (157.5:2)node{$\bullet$} to (202.5:2)node{$\bullet$}
to (247.5:2)node{$\bullet$} to (292.5:2)node{$\bullet$} to (-22.5:2)node{$\bullet$} to (22.5:2);
\draw[thick] (22.5:2)node[right]{4} (67.5:2)node[above]{5} (112.5:2)node[above]{6} (157.5:2)node[left]{7} (202.5:2)node[left]{8}
(247.5:2)node[below]{1} (292.5:2)node[below]{2} (-22.5:2)node[right]{3};
\draw[thick] (202.5:2) to (292.5:2) to (157.5:2)  (-22.5:2) to (202.5:2);
\draw[thick] (112.5:2) to (22.5:2) to (157.5:2) to (67.5:2);
\draw[thick, red] (157.5:2) to (-22.5:2);
\end{tikzpicture}
\qquad
\begin{tikzpicture}[scale=.8]
\draw[thick] (0,0)node{$\xrightarrow{\text{$\{\{3,7\}\}$-rotation}}$} (-90:2.5)node{ };
\end{tikzpicture}
\qquad
\begin{tikzpicture}[scale=.8]
\draw[thick] (22.5:2)node{$\bullet$} to (67.5:2)node{$\bullet$} to (112.5:2)node{$\bullet$}  to (157.5:2)node{$\bullet$} to (202.5:2)node{$\bullet$}
to (247.5:2)node{$\bullet$} to (292.5:2)node{$\bullet$} to (-22.5:2)node{$\bullet$} to (22.5:2);
\draw[thick] (22.5:2)node[right]{4} (67.5:2)node[above]{5} (112.5:2)node[above]{6} (157.5:2)node[left]{7} (202.5:2)node[left]{8}
(247.5:2)node[below]{1} (292.5:2)node[below]{2} (-22.5:2)node[right]{3};
\draw[thick] (-22.5:2) to (247.5:2) to (157.5:2) to (292.5:2);
\draw[thick] (22.5:2) to (112.5:2) to (-22.5:2) to (67.5:2);
\draw[thick, red] (157.5:2) to (-22.5:2);
\end{tikzpicture}
\caption{An example of rotation of Ptolemy diagrams of a polygon}
\label{fig:rP}
\end{figure}
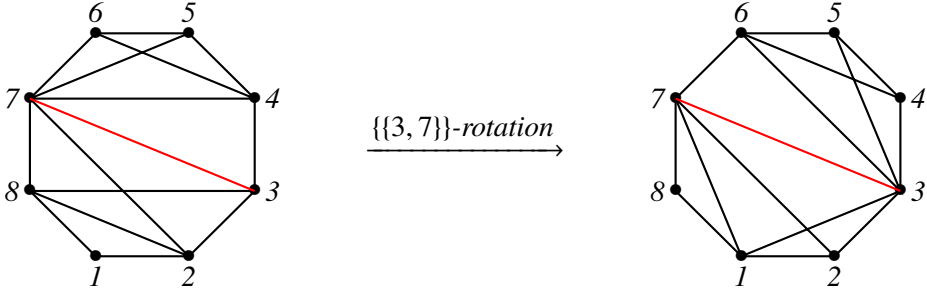
\end{exm}

The following lemma is essentially from \cite{CCS}.

\begin{lem}\label{lem:b1}
Let $\{i,j\}$ and $\{k,l\}$ be diagonals of $P_{n+3}$. Then they cross each other if and only if $\Hom_\CA(X_{k,l},X_{i,j}[1])\neq 0$. Moreover, in this case, there is a non-split triangle
\begin{equation}\label{eq:b1}
X_{i,j}\rightarrow X_{k,j}\oplus X_{i,l}\rightarrow X_{k,l}\rightarrow X_{i,j}[1]
\end{equation}
where $X_{a,b}$ is taken to be a zero object if $\{a,b\}$ is an edge of $P_{n+3}$.
\end{lem}

Now we prove that rotation of Ptolemy diagrams gives a geometric model of mutation of torsion pairs.

\begin{thm}\label{thm:An}
Let $\UU$ be a Ptolemy diagram of $P_{n+3}$. Then $X(F_\UU)=I(X(\UU))$ and the map $X$ induces a bijection between subsets of $F_\UU$ and functorially finite subcategories of $\C_{A_n}$ contained in $I(X(\UU))$. Moreover, for any subset $\DD\subset F_\UU$, we have that \[\mu^{-1}(X(\UU);X(\DD))=X\left(\rho(\UU;\DD)\right).\] Consequently, a rotation of a Ptolemy diagram is again a Ptolemy diagram.
\end{thm}

\begin{proof}
Since $\C_{A_n}$ is Hom-finite and contains finitely many indecomposable objects up to isomorphism, each of its subcategories is functorially finite. Then the first assertion follows directly from Lemma~\ref{lem:b1}.

For the second assertion, it is sufficient to prove that for any diagonal $\{i,j\}$ in $\UU$ which is not in $\DD$, we have $\mu^{-1}_{X(\DD)}(X_{i,j})=X_{k,l}$ where $\{k,l\}=\rho_\DD(\{i,j\}).$ Let $P'$ be the $\DD$-cell containing $\{i,j\}$. By definition, we have that $\{k,l\}$ is the rotation of $\{i,j\}$ in $P'$. So $\{k,l\}$ is a diagonal of $P'$ crossing $\{i,j\}$, and all of $\{i,k\},\{k,j\},\{j,l\},\{l,i\}$ are edges of $P'$. Then all of $X_{i,k},X_{k,j},X_{j,l},X_{l,i}$ belong to $X(\DD)$ and by Lemma~\ref{lem:b1} we have $X_{k,l}\in {}^\bot X(\DD)[1]$. Applying $\Hom_\CA(-,D)$, for $D\in X(\DD)$, to the triangle \eqref{eq:b1} in Lemma~\ref{lem:b1}, we have an exact sequence
\[ \Hom_\CA(X_{k,j}\oplus X_{i,l},D)\rightarrow \Hom_\CA(X_{i,j},D)\rightarrow \Hom_\CA(X_{k,l},D),\]
where the last term is zero by$X_{k,l}\in {}^\bot X(\DD)[1]$. This implies that the first morphism in \eqref{eq:b1} is a left $X(\DD)$-approximation and it is left minimal since $X_{k,l}$ is indecomposable. So by definition, $X_{k,l}=\mu^{-1}_{X(\DD)}(X_{i,j})$ and we are done. Then the last assertion follows from Proposition~\ref{prop:hjr}.
\end{proof}

\subsection{Infinite Dynkin type A}

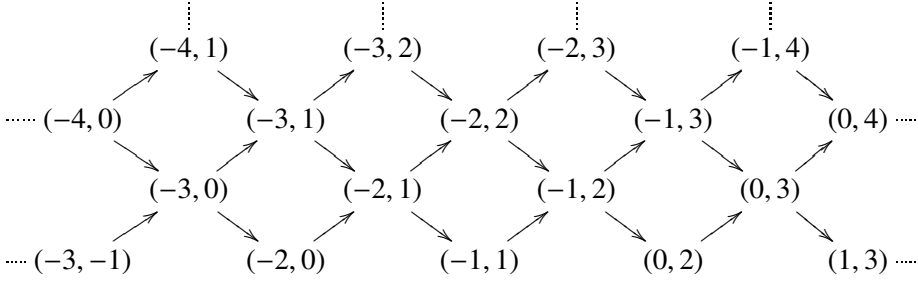
\begin{figure}[ht]\centering
\[\xymatrix@C=0em@R=1em{
&&&&&&&&&&&&\\
&&&(-4,1)\ar[rd]\ar@{..}[u]&&(-3,2)\ar[rd]\ar@{..}[u]&&(-2,3)\ar[rd]\ar@{..}[u]&&(-1,4)\ar[rd]\ar@{..}[u]&&&\\
&&(-4,0)\ar@{..}[ll]\ar[rd]\ar[ru]&&(-3,1)\ar[rd]\ar[ru]&&(-2,2)\ar[rd]\ar[ru]&&(-1,3)\ar[rd]\ar[ru]&&(0,4)\ar@{..}[rr]&&\\
&&&(-3,0)\ar[rd]\ar[ru]&&(-2,1)\ar[rd]\ar[ru]&&(-1,2)\ar[rd]\ar[ru]&&(0,3\ar[rd]\ar[ru])&&&\\
&&(-3,-1)\ar@{..}[ll]\ar[ru]&&(-2,0)\ar[ru]&&(-1,1)\ar[ru]&&(0,2)\ar[ru]&&(1,3)\ar@{..}[rr]&&
}\]
\caption{The Auslander-Reiten quiver of $C_{A_\infty}$}
\label{fig:AR}
\end{figure}

The cluster category $\C_{A_{\infty}}$ of type $A_\infty$ was studied in \cite{HJ} and \cite{Ng}. The Auslander-Reiten quiver of $\C_{A_\infty}$ is $\mathbb{Z}A_\infty$, see Figure~\ref{fig:AR}, where the (isoclasses of) indecomposable objects are indexed by the pairs $(i,j)$ of integers with $j-i\geq2$. The action of shift functor [1] on indecomposable objects is given by $(i,j)[1]=(i-1,j-1)$.

Let $P_\infty$ be an infinity-gon whose vertices labeled by the integers. Similarly as in a (finite) polygon, we use $\{i,j\}$ to denote the diagonal or edge of $P_\infty$ linking $i$ and $j$, and for any diagonal $\{i,j\}$, its rotation in $P_\infty$ is defined to be $\{i-1,j-1\}$. Then the map $\{i,j\}\mapsto X_{i,j}:=(\min\{i,j\},\max\{i,j\})$ gives a bijection from the set of diagonals of $P_\infty$ to the set $\ind\C_{A_\infty}$. This bijection induces a bijection between sets of diagonals of $P_\infty$ and subcategories of $\C_{A_\infty}$. We still use $X(\UU)$ to denote the subcategory associated to a set $\UU$ of diagonals.

We recall some definitions and results from \cite{HJ,Ng}.

\begin{defn}
Let $\UU$ be a set of diagonals of $P_\infty$.
\begin{enumerate}
  \item[(1)] An integer $n$ is called a left (resp. right) fountain of $\UU$ if $\UU$ contains infinitely many diagonals of the form $\{m,n\}$ for $m<n$ (resp. for $m>n$).
  \item[(2)] An integer $n$ is called a fountain of $\UU$, if it is both a left and a right fountain of $\UU$.
  \item[(3)] The set $\UU$ is called locally finite, if for each integer $n$, $n$ is neither a left nor a right fountain of $\UU$.
  \item[(4)] The set $\UU$ is called a Ptolemy diagram, if the following three conditions hold:
  \begin{enumerate}
    \item[(i)] each right fountain of $\mathfrak{U}$ is a fountain;
    \item[(ii)] if $i$ is a left fountain of $\UU$, $j$ is a right fountain of $\UU$ and $j-i\geq2$, then $\{i,j\}\in\UU$.
    \item[(iii)] for each pair of crossing diagonals $\{i,j\}$ and $\{k,l\}$ in $\mathfrak{U}$, those of the pairs $\{i,k\}$, $\{k,j\}$, $\{j,l\}$ and $\{l,i\}$ which are diagonals belong to $\mathfrak{U}$. 
  \end{enumerate}
\end{enumerate}
\end{defn}

\begin{prop}\label{prop:basic}
The following hold.
\begin{enumerate}
  \item[(1)] Two diagonals $(i,j)$ and $(k,l)$ cross each other if and only if $\Hom_{\C_{A_\infty}}(X_{i,j},X_{k,l}[1])\neq0$.
  \item[(2)] A set $\UU$ of diagonals is locally finite or has a fountain if and only if $X(\UU)$ is functorially finite.
  \item[(3)] A set $\UU$ of diagonals is a Ptolemy diagram if and only if $(X(\UU),X(\UU)^\bot)$ is a torsion pair in $\C_{A_\infty}$.
\end{enumerate}
\end{prop}

\begin{proof}
The first assertion follows from \cite[Lemma~3.6, Corollary~2.3]{HJ}. The second assertion follows from \cite[Theorem~2.2, Theorem~2.3]{Ng} and the third one is exactly \cite[Theorem~3.18]{Ng}.
\end{proof}

In the following, we define frame and mutation of a Ptolemy diagram of an infinity-gon.

\begin{defn}
Let $\UU$ be a Ptolemy diagram. The set of diagonals in $\UU$ that does not cross any diagonals in $\UU$ is called the frame of $\UU$, denoted by $F_\UU$.
\end{defn}

\begin{defn}
Let $\DD$ be a set of non-crossing diagonals such that it is locally finite or has a fountain.
\begin{itemize}
  \item[(1)] The diagonals in $\DD$ divide the polygon $P_\infty$ into (maybe infinitely) many polygons and in addition an infinity-gon when $\DD$ is locally finite, called $\DD$-cells. Then any diagonal $\{i,j\}$ which neither is in $\DD$ nor crosses any diagonals in $\DD$ is a diagonal of a $\DD$-cell. We call the rotation of $\{i,j\}$ in this $\DD$-cell the $\DD$-rotation of $\{i,j\}$ 
  and denote it by $\rho_\D(\{i,j\})$.
  \item[(2)] Let $\UU$ be a Ptolemy diagram with $\DD\subset F_\UU$. The $\DD$-rotation of $\UU$ is defined as
  \[\rho(\UU;\DD):=\{\rho_\DD(\{i,j\})\mid \{i,j\}\in\UU\setminus \DD\}\cup \DD.\]
\end{itemize}
\end{defn}

\begin{exm}
Let $\UU=\{\{1,3\},\{1,4\},\{2,4\},\{1,5\}\}\cup\{\{t,1\}\mid t\leq-1\}$, see the upper picture in Figure~\ref{fig:rP2}. Then $\UU$ is a Ptolemy diagram with $F_\UU=\{\{1,4\},\{1,5\}\}\cup\{\{t,1\}\mid t\leq-1\}$. Note that $F_\UU$ neither is locally finite nor has a fountain.
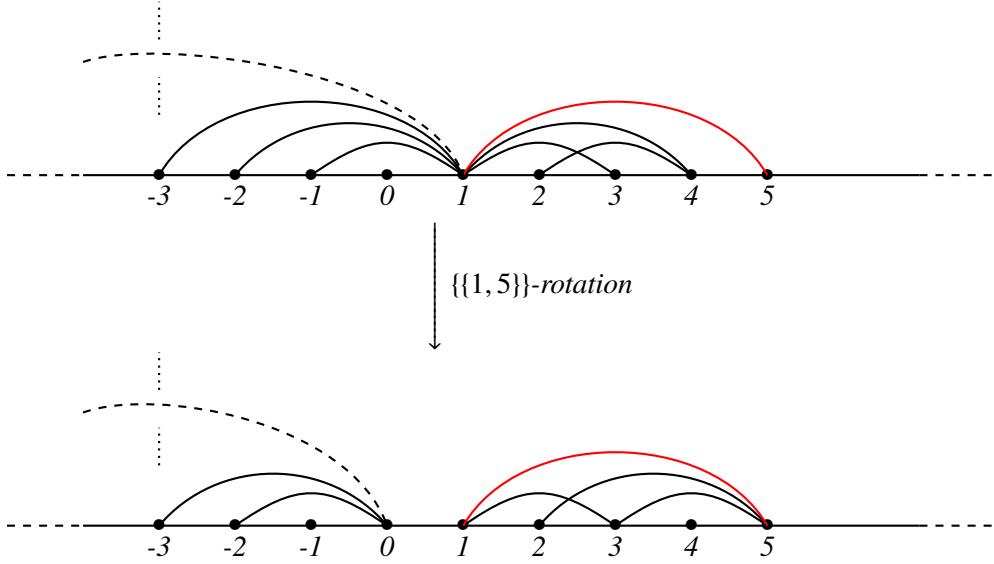
\begin{figure}[ht]\centering
\begin{tikzpicture}
\draw[thick] (-4,0) to (7,0);
\draw[thick, dashed] (-5,0) to (-4,0) (7,0) to (8,0);
\draw[thick] (-3,0)node{$\bullet$}node[below]{-3} (-2,0)node{$\bullet$}node[below]{-2} (-1,0)node{$\bullet$}node[below]{-1} (0,0)node{$\bullet$}node[below]{0}
(1,0)node{$\bullet$}node[below]{1} (2,0)node{$\bullet$}node[below]{2} (3,0)node{$\bullet$}node[below]{3} (4,0)node{$\bullet$}node[below]{4} (5,0)node{$\bullet$}node[below]{5};
\draw[thick] (1,0) ..  controls +(35:1) and +(145:1) .. (3,0);
\draw[thick] (1,0) ..  controls +(50:1.2) and +(130:1.2) .. (4,0);
\draw[thick] (2,0) ..  controls +(35:1) and +(145:1) .. (4,0);
\draw[thick, red] (1,0) ..  controls +(60:1.5) and +(120:1.5) .. (5,0);
\draw[thick] (-1,0) ..  controls +(35:1) and +(145:1) .. (1,0);
\draw[thick] (-2,0) ..  controls +(50:1.2) and +(130:1.2) .. (1,0);
\draw[thick] (-3,0) ..  controls +(60:1.5) and +(120:1.5) .. (1,0);
\draw[thick, dashed] (-4,1.5) ..  controls +(20:1) and +(110:1.5) .. (1,0);
\draw[thick, dotted] (-3,.8) to (-3,1.3) (-3,1.8) to (-3,2.3);
\end{tikzpicture}

$\qquad\left\downarrow\rule{0cm}{1cm}\right.\{\{1,5\}\}\text{-rotation}$

\begin{tikzpicture}
\draw[thick] (-4,0) to (7,0);
\draw[thick, dashed] (-5,0) to (-4,0) (7,0) to (8,0);
\draw[thick] (-3,0)node{$\bullet$}node[below]{-3} (-2,0)node{$\bullet$}node[below]{-2} (-1,0)node{$\bullet$}node[below]{-1} (0,0)node{$\bullet$}node[below]{0}
(1,0)node{$\bullet$}node[below]{1} (2,0)node{$\bullet$}node[below]{2} (3,0)node{$\bullet$}node[below]{3} (4,0)node{$\bullet$}node[below]{4} (5,0)node{$\bullet$}node[below]{5};
\draw[thick] (1,0) ..  controls +(35:1) and +(145:1) .. (3,0);
\draw[thick] (2,0) ..  controls +(50:1.2) and +(130:1.2) .. (5,0);
\draw[thick] (3,0) ..  controls +(35:1) and +(145:1) .. (5,0);
\draw[thick, red] (1,0) ..  controls +(60:1.5) and +(120:1.5) .. (5,0);
\draw[thick] (-2,0) ..  controls +(35:1) and +(145:1) .. (0,0);
\draw[thick] (-3,0) ..  controls +(50:1.2) and +(130:1.2) .. (0,0);
\draw[thick, dashed] (-4,1.5) ..  controls +(20:1) and +(110:1.5) .. (0,0);
\draw[thick, dotted] (-3,.8) to (-3,1.3) (-3,1.8) to (-3,2.3);
\end{tikzpicture}
\caption{An example of rotation of Ptolemy diagrams of $P_\infty$}
\label{fig:rP2}
\end{figure}
Let $\DD=\{\{1,5\}\}\subset F_\UU$ which is locally finite. Then $\DD$ divide $P_\infty$ into a polygon whose vertices are labeled by $1, 2, 3, 4, 5$, and an infinity-gon whose vertices are labeled by $\cdots, -3,-2,-1,0,1,5,6,\cdots$. Then the $\DD$-rotation of $\UU$ is the lower picture in Figure~\ref{fig:rP2}, where for example, $\rho_\DD(\{-2,1\})=\{-3,0\}$ and $\rho_\DD(\{1,4\})=\{5,2\}$. It is clear that $\rho(\UU;\DD)$ is again a Ptolemy diagram.
\end{exm}

\begin{lem}\label{lem:tri}
Let $a,b,c,d$ be four integers with $a<c<b<d$. Then there are two non-split triangles in $\C_{A_\infty}$:
\[(a,b)\rightarrow (a,d)\oplus(c,b)\rightarrow(c,d)\rightarrow(a,b)[1]\]
and
\[(c,d)\rightarrow (a,c)\oplus(b,d)\rightarrow(a,b)\rightarrow(c,d)[1],\]
where $(x,y)$ is taken to be a zero object when $y-x=1$.
\end{lem}

\begin{proof}
The first non-split triangle follows straightly from the Auslander-Reiten quiver of $\C_{A_\infty}$. Now we show the existence of the second one. Consider the following two non-split triangles in $\C_{A_\infty}$ from the Auslander-Reiten quiver of $\C_{A_\infty}$:
\[(c,b+1)\rightarrow(c,d)\rightarrow(b,d)\rightarrow(c-1,b),\]
and
\[(c,b+1)\rightarrow(a,c)\rightarrow(a,b)\rightarrow(c-1,b).\]
Applying $\Hom((b,d),-)$ to the last triangle, we have the following exact sequence
\[\Hom((b,d),(a,b))\rightarrow\Hom((b,d),(c-1,b))\rightarrow\Hom((b,d),(a,c)[1]).\]
where the third term is zero by Proposition~\ref{prop:basic} since $(a,c)$ does not cross $(b,d)$. So the first map in this exact sequence is surjective, which means that any morphism from $(b,d)$ to $(c-1,b)$ factors through $(a,b)$. Then there is a communicative diagram
$$\begin{array}{ccccccc}
(c,b+1)&\rightarrow&(c,d)&\rightarrow&(b,d)&\rightarrow&(c-1,b)\\
\parallel&&\downarrow &&\downarrow&&\parallel\\
(c,b+1)&\rightarrow&(a,c)&\rightarrow&(a,b)&\rightarrow&(c-1,b).
\end{array}$$
By \cite[Lemma 2.2]{XZ1}, we get the second triangle in the lemma. It is non-split since $(c,d)\oplus (a,b)\ncong (a,c)\oplus
(b,d)$.
\end{proof}

We have the following analogue result for $\C_{A_\infty}$. 

\begin{thm}\label{thm:Ainf}
Let $\UU$ be a Ptolemy diagram of $P_\infty$. Then $X(F_\UU)=I(X(\UU))$ and the map $X$ induces a bijection from subsets of $F_\UU$ which are locally finite or have a fountain to functorially finite subcategories of $\C_{A_\infty}$ contained in $I(X(\UU))$. Moreover, for any subset $\DD\subset F_\UU$ which is locally finite or have a fountain, we have that 
\[\mu^{-1}(X(\UU);X(\DD))=X\left( \rho(\UU;\DD)\right).\] 
Consequently, a rotation of a Ptolemy diagram is again a Ptolemy diagram.
\end{thm}

\begin{proof}
The first assertion follows straightly from Proposition~\ref{prop:basic} (1) and (2). The second assertion can be proved by a similar way with the proof of Theorem~\ref{thm:An}, using Lemma~\ref{lem:tri}. Then the last assertion follows from Proposition~\ref{prop:basic}(3).

\end{proof}

\end{document}